\documentclass[a4paper,12pt]{article}
\usepackage[utf8]{inputenc}
\usepackage[UKenglish]{isodate}
\cleanlookdateon
\usepackage{amsthm}
\usepackage{amssymb}
\usepackage{amsmath}
\usepackage{enumitem}
\usepackage[T1]{fontenc}
\usepackage{xcolor}
\usepackage{graphicx}
\graphicspath{ {./images/} }
\usepackage{caption}
\usepackage{subcaption}
\usepackage{authblk}
\usepackage{hyperref}

\hypersetup{
  colorlinks   = true,
  urlcolor     = blue!50!black,
  linkcolor    = blue!50!black,
  citecolor   = blue!50!black
}

\newcommand{\floor}[1]{\left \lfloor{ #1 }\right \rfloor}

\title{\vspace{-0.67cm}Planar graphs with the maximum number of induced 6-cycles}
\date{22 January 2024}
\author{Michael Savery}

\affil{Mathematical Institute, University of Oxford, Oxford OX2 6GG, UK, and the Heilbronn Institute for Mathematical Research, Bristol, UK\\ \texttt{savery@maths.ox.ac.uk}}

\newtheorem{theorem}{Theorem}
\newtheorem{lemma}[theorem]{Lemma}
\newtheorem{corollary}[theorem]{Corollary}
\newtheorem{claim}{Claim}
\newtheorem{conjecture}[theorem]{Conjecture}

\theoremstyle{definition}
\newtheorem{definition}[theorem]{Definition}
\newtheorem{fact}[theorem]{Fact}

\advance\textwidth2\evensidemargin
\begin{document}

\maketitle

\begin{abstract}
    \noindent For large $n$ we determine the maximum number of induced 6-cycles which can be contained in a planar graph on $n$ vertices, and we classify the graphs which achieve this maximum. In particular we show that the maximum is achieved by the graph obtained by blowing up three pairwise non-adjacent vertices in a 6-cycle to sets of as even size as possible, and that every extremal example closely resembles this graph. This extends previous work by the author which solves the problem for 4-cycles and 5-cycles. The 5-cycle problem was also solved independently by Ghosh, Gy\H{o}ri, Janzer, Paulos, Salia, and Zamora.
\end{abstract}

\section{Introduction}\label{sec:intro}

The problem of determining for a fixed small graph $H$ the maximum number of induced copies of $H$ which can be contained in a graph on $n$ vertices was first considered in 1975 by Pippenger and Golumbic~\cite{pippenger}, and has since received considerable attention. As summarised in~\cite{even-zohal}, this maximum is now known asymptotically for all graphs $H$ on at most four vertices except the path of length 3. In the case where $H$ is a $k$-cycle for $k\geq 5$, Pippenger and Golumbic conjectured in their paper that the maximum is asymptotically $\frac{n^k}{k^k-k}$, which they showed is an asymptotic lower bound for all $k$-vertex graphs $H$. In 2016 this conjecture was verified in the case $k=5$ by Balogh, Hu, Lidick\'{y}, and Pfender in~\cite{balogh}. The current best known upper bound for general $k$-cycles is $2n^k/k^k$, which is due to Kr\'{a}l', Norin, and Volec~\cite{kral}. Results for other graphs $H$ can be found in~\cite{hatami} and~\cite{yuster}.

A related problem on which significant progress has been made recently is to consider what happens when we restrict to planar graphs. There are two interesting problems in this setting, firstly to determine the maximum number of not necessarily induced copies of a small graph $H$ that can be contained in a planar graph on $n$ vertices, which we write as $f(n,H)$, and secondly to determine the same quantity when we insist that the copies of $H$ are induced. We write $f_I(n,H)$ for this second quantity.

In a moment we shall focus on the case where $H$ is a cycle, but we first note that $f(n,H)$ has also been studied for other graphs $H$ (see, for example,~\cite{alon},~\cite{cox2021counting},~\cite{eppstein},~\cite{ghoshp5},~\cite{grzesik}, \cite{wood}, and~\cite{wormald}). In the general case, Huynh, Joret, and Wood determined in~\cite{huynh2021subgraph} the order of magnitude of $f(n,H)$ for all $H$ in terms of a graph parameter called the `flap-number' of $H$. Liu then showed in~\cite{liu2021homomorphism} (Corollary 6.1 in that paper) that for all graphs~$H$, $f_I(n,H)=\Theta(f(n,H))$, which in turn determines the order of magnitude of $f_I(n,H)$ for all $H$.

Focussing now on the case where $H$ is a cycle, in 1979 Hakimi and Schmeichel~\cite{hakimi} showed that $f(n,C_3)=3n-8$ for $n\geq 3$ and $f(n,C_4)=\frac{1}{2}(n^2+3n-22)$ for $n\geq 4$. Much more recently, in 2019, $f(n,C_5)$ was determined for all $n\geq 5$ by Gy\H{o}ri, Paulos, Salia, Tompkins, and Zamora in~\cite{gyori}. Since every copy of $C_3$ in a graph is an induced copy, Hakimi and Schmeichel's result also determines $f_I(n,C_3)$ for all $n\geq 3$. The values of $f_I(n,C_4)$ and $f_I(n,C_5)$ were determined for large enough $n$ in~\cite{savery2021planar}. In each case the extremal graphs were classified in the respective paper, with a small correction to the extremal graphs for $f(7,C_4)$ and $f(8,C_4)$ given by Alameddine in~\cite{alameddine}. The value of $f_I(n,C_5)$ for large $n$ was also found independently by Ghosh, Gy\H{o}ri, Janzer, Paulos, Salia, and Zamora in~\cite{ghosh2021maximum}.

\begin{figure}
    \centering
    \begin{subfigure}[h]{0.45\textwidth}
        \centering
        \includegraphics[width=\textwidth]{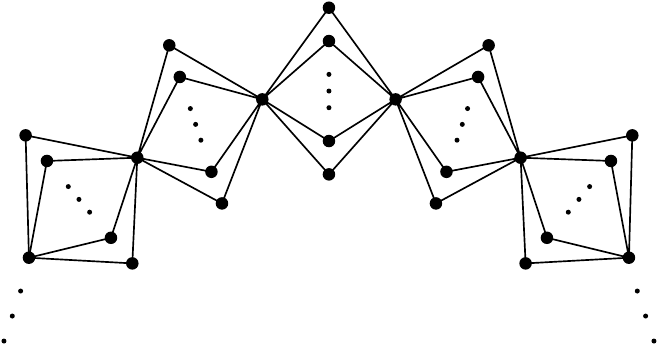}
        \caption{Graph $F_{n,m}$ when $m$ is even}
        \label{fig:F_nm_even}
     \end{subfigure}
     \hfill
     \begin{subfigure}[h]{0.45\textwidth}
        \centering
        \includegraphics[width=\textwidth]{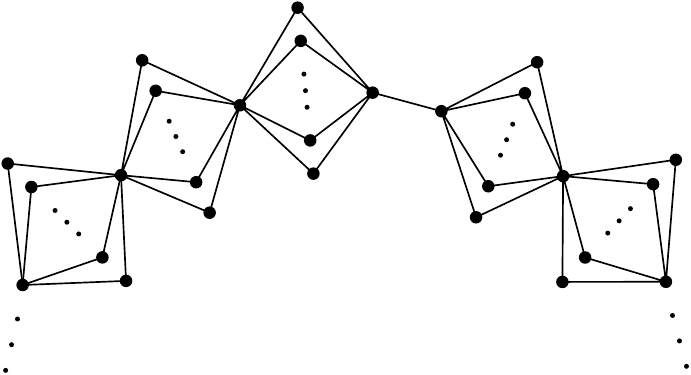}
        \caption{Graph $F_{n,m}$ when $m$ is odd}
        \label{fig:F_nm_odd}
     \end{subfigure}
     
     \begin{subfigure}[h]{0.45\textwidth}
        \centering
        \includegraphics[width=\textwidth]{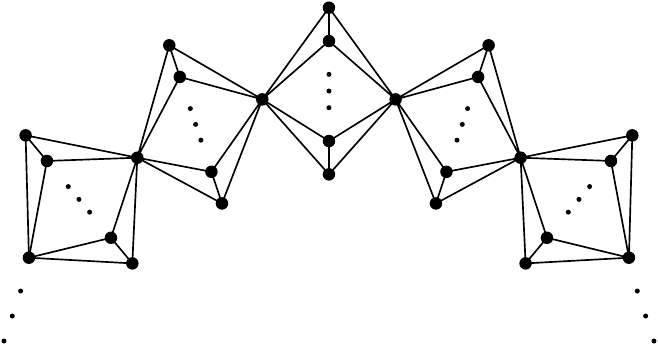}
        \caption{Graph $F'_{n,m}$ when $m$ is even}
        \label{fig:Fprime_nm_even}
     \end{subfigure}
     \hfill
     \begin{subfigure}[h]{0.45\textwidth}
        \centering
        \includegraphics[width=\textwidth]{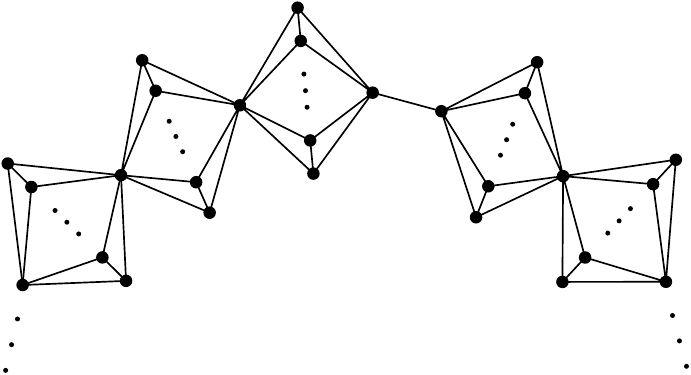}
        \caption{Graph $F'_{n,m}$ when $m$ is odd}
        \label{fig:Fprime_nm_odd}
     \end{subfigure}
        \caption{The graphs $F_{n,m}$ and $F'_{n,m}$ in different cases}
        \label{fig:Fnm}
\end{figure}

No exact results are known for longer cycles, but a recent string of papers have yielded impressive \emph{asymptotic} results. This progress was initiated by Cox and Martin~\cite{cox2021counting}, who developed a method for upper bounding $f(n,H)$ for even cycles $H$ via a reduction to a maximum likelihood estimator problem on graphs. They went on to solve this problem for small cycles~\cite{cox2021counting,cox2021maximum}, before Lv, Gy\H{o}ri, He, Salia, Tompkins, and Zhu~\cite{lv2022maximum} solved it in general to show that $f(n,C_{2k})=\left(\frac{n}{k}\right)^k+o(n^k)$ for all $k\geq 3$. Planar graphs achieving these asymptotics were first exhibited in~\cite{gyori2020generalized}: for $m\geq 3$, let $F_{n,m}$ be the $n$-vertex graph obtained by blowing up $\floor{\frac{m}{2}}$ pairwise non-adjacent vertices in an $m$-cycle to sets of as even size as possible (see Figures~\hyperref[fig:F_nm_even]{\ref*{fig:Fnm}(\subref*{fig:F_nm_even})} and~\hyperref[fig:F_nm_odd]{\ref*{fig:Fnm}(\subref*{fig:F_nm_odd})}), then $F_{n,2k}$ is a planar $n$-vertex graph containing $\left(\frac{n}{k}\right)^k+o(n^k)$ cycles of length $2k$. In fact, all of these $2k$-cycles are induced, so since $f_I(n,H)\leq f(n,H)$ for all $n$ and $H$, we have also that $f_I(n,C_{2k})=\left(\frac{n}{k}\right)^k+o(n^k)$ for all $k\geq 3$.

In the case where $H$ is an odd cycle of length greater than 5, we can obtain a lower bound on $f_I(n,C_{2k+1})$ in the same spirit as above by observing that $F_{n,2k+1}$ is an $n$-vertex planar graph containing $\left(\frac{n}{k}\right)^k+o(n^k)$ induced $(2k+1)$-cycles. To obtain a better lower bound on $f(n,C_{2k+1})$ we can consider the following graphs, first defined in~\cite{ghoshp5}: for each $m\geq 3$, let $F'_{n,m}$ be the graph obtained from $F_{n,m}$ by adding a path through each of the blown-up vertex classes (see Figures~\hyperref[fig:Fprime_nm_even]{\ref*{fig:Fnm}(\subref*{fig:Fprime_nm_even})} and~\hyperref[fig:Fprime_nm_odd]{\ref*{fig:Fnm}(\subref*{fig:Fprime_nm_odd})}). Then $F'_{n,2k}$ is an $n$-vertex planar graph containing $2k\left(\frac{n}{k}\right)^k+o(n^k)$ (non-induced) $(2k+1)$-cycles. Very recently, Heath, Martin, and Wells~\cite{heath23maximum} extended Cox and Martin's reduction lemmas to the case of odd cycles, and showed that for $k\in\{3,4\}$ the graphs $F'_{n,2k}$ are asymptotically optimal for the non-induced problem, that is, $f(n,C_{2k+1})=2k\left(\frac{n}{k}\right)^k+o(n^k)$ for $k\in\{3,4\}$. They conjectured that this extends to all $k\geq 5$, and showed a general upper bound of $f(n,C_{2k+1})\leq 3k\left(\frac{n}{k}\right)^k+o(n^k)$. No better upper bounds on $f_I(n,C_{2k+1})$ are known for~$k\geq 3$. \pagebreak

In~\cite{savery2021planar} the following conjecture was made regarding the exact solution of the problem for induced even cycles of length at least $6$.

\begin{conjecture}[\cite{savery2021planar}]\label{conj:induced_even}
For $k\geq 6$ and $n$ sufficiently large relative to $k$, the graph $F_{n,k}$ contains $f_I(n,C_k)$ induced $k$-cycles.
\end{conjecture}

In this paper we will show that for large enough $n$, the $n$-vertex planar graphs containing $f_I(n,C_6)$ induced 6-cycles are exactly those which are subgraphs of $F'_{n,6}$ and which contain $F_{n,6}$ as a subgraph, thus establishing the conjecture when $k=6$. By counting the number of induced 6-cycles in $F_{n,6}$ we can obtain from this a closed-form expression for $f_I(n,C_6)$ for large $n$. Our approach is different from that introduced by Cox and Martin, which in its current form is only of use for obtaining asymptotic results.

Before stating our result, for each $n\geq 6$ we will define $\mathcal{F}_n$ to be the family consisting of the subgraphs of $F'_{n,6}$ which contain $F_{n,6}$ as a subgraph. Formally, we define $\mathcal{F}_n$ as follows.

\begin{definition}\label{def:Fn}
For $n\geq 6$ the family $\mathcal{F}_n$ consists of the graphs $G$ such that all of the following hold:
\begin{enumerate}[label=(\roman*)]
    \item $G$ has exactly $n$ vertices, with vertex set $\{u_1,u_2,u_3\}\cup A \cup B\cup C$, where $|A|+|B|+|C|=n-3$ and the sizes of $A$, $B$, and $C$ are as equal as possible,
    \item for all $a\in A$, $b\in B$, and $c\in C$, $G$ contains the edges $au_1$, $au_2$, $bu_2$, $bu_3$, $cu_3$, and $cu_1$, and
    \item there exist labellings of $A$, $B$, and $C$ as, respectively, $\{a_1,\dots,a_{|A|}\}$, $\{b_1,\dots,b_{|B|}\}$, and $\{c_1,\dots,c_{|C|}\}$ such that any remaining edges of $G$ form a subset of $\{a_ia_{i+1}:1\leq i\leq |A|-1\}\cup \{b_ib_{i+1}:1\leq i\leq |B|-1\}\cup \{c_ic_{i+1}:1\leq i\leq |C|-1\}$.
\end{enumerate}
\end{definition}

\begin{figure}[!ht]
    \centering
    \includegraphics[width=0.5\textwidth]{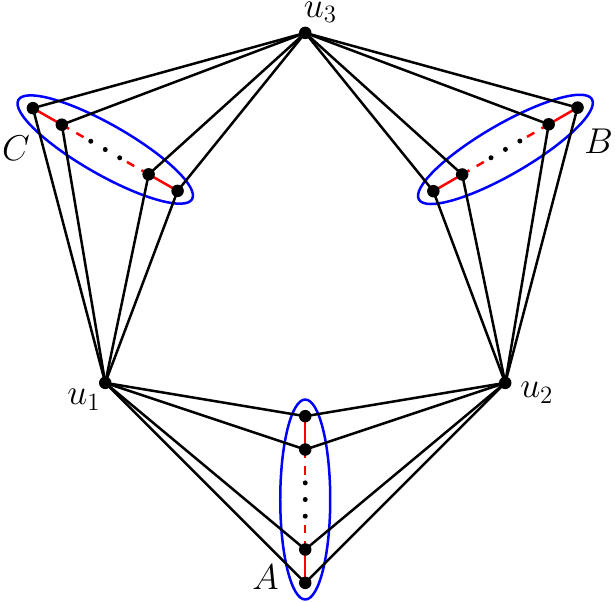}
    \caption{A graph in $\mathcal{F}_n$}
    \label{fig:Fn}
\end{figure}

The structure of a graph in $\mathcal{F}_n$ is illustrated in Figure~\ref{fig:Fn}. In this illustration, red lines indicate edges each of whose addition or omission does not affect whether the graph is in $\mathcal{F}_n$. We are now ready to state the main result of the paper.

\begin{theorem}\label{thm:main}
For large $n$, $G$ is an $n$-vertex planar graph containing $f_I(n,C_6)$ induced 6-cycles if, and only if, $G\in \mathcal{F}_n$. Hence for large $n$, 
\[
f_I(n,C_6)=
\begin{cases}
\left(\frac{n}{3}-1\right)^3=\frac{1}{27}(n^3-9n^2+27n-27), & \text{if $n\equiv 0\pmod{3}$} \\
\left(\frac{n-4}{3}\right)^2\left(\frac{n-1}{3}\right)=\frac{1}{27}(n^3-9n^2+24n-16), & \text{if $n\equiv 1\pmod{3}$} \\ 
\left(\frac{n-2}{3}\right)^2\left(\frac{n-5}{3}\right)=\frac{1}{27}(n^3-9n^2+24n-20), & \text{if $n\equiv 2\pmod{3}$}.
\end{cases}
\] 
\end{theorem}

\subsection{Notation and outline of the paper}\label{sec:outline}
We use standard notation from graph theory in this paper. In particular, $C_k$ is the $k$-cycle graph, $K_{a,b}$ is the complete bipartite graph with parts of size $a$ and $b$, and for a graph~$G$ and a subset $S$ of its vertex set, $G-S$ is the subgraph of $G$ obtained by deleting the vertices in $S$. For a vertex $v$ of a graph $G$ we write $N(v)$ for the (open) neighbourhood of $v$ in $G$, that is, the set of vertices of $G$ which are adjacent to $v$.

In the next section we first introduce some of the framework we will use to prove Theorem~\ref{thm:main}, and then we find some structure in large planar graphs in which no vertex is in `few' induced 6-cycles. In Section~\ref{sec:vtx_in_few} we use this to show that every large planar graph contains a vertex in `not too many' induced 6-cycles, and that if every vertex is in `many' induced 6-cycles and the number of vertices in the graph is a multiple of 3, then the graph is in $\mathcal{F}_n$ for the appropriate $n$. In Section~\ref{sec:main_proof} we use these results to show that for large $n$ every vertex in a planar $n$-vertex graph containing $f_I(n,C_6)$ induced 6-cycles is in `many' induced 6-cycles, from which we can deduce that all such graphs are in $\mathcal{F}_n$, which is the substance of Theorem~\ref{thm:main}. Finally in Section~\ref{sec:longer_cycles} we consider to what extent the arguments used are applicable to the problem of determining $f_I(n,C_k)$ for $k>6$.

\section{Preliminaries}\label{sec:prelims}

We will use the following notation throughout this paper. Let $v$ be a vertex of a planar graph $G$ with distinct neighbours $u$ and $w$. Define $X^{uvw}$ to be the set of vertices of $G-\{u,v,w\}$ which are in an induced 6-cycle in $G$ containing the path $uvw$. Then let $X_1^{uvw}=X^{uvw}\cap N(u)$, $X_3^{uvw}=X^{uvw}\cap N(w)$, and $X_2^{uvw}=X^{uvw}\setminus(X_1^{uvw}\cup X_3^{uvw})$, so that every induced 6-cycle in $G$ containing the path $uvw$ is of the form $uvwx_3x_2x_1$ for some $x_1\in X_1$, $x_2\in X_2$, and $x_3\in X_3$.

Clearly the sets $X_1^{uvw}$, $X_2^{uvw}$, and $X_3^{uvw}$ are pairwise disjoint. If there exists an induced 6-cycle in $G$ containing the path $uvw$, then they are all non-empty and every vertex in $X_1^{uvw}$ has a neighbour in $X_2^{uvw}$, every vertex in $X_2^{uvw}$ has a neighbour in $X_1^{uvw}$ and a neighbour in $X_3^{uvw}$, and every vertex in $X_3^{uvw}$ has a neighbour in $X_2^{uvw}$. The following fundamental lemma gives further properties of these sets and introduces a framework that we will use often. It is adapted from Lemma~1 in~\cite{ghosh2021maximum}, which is a similar result for 5-cycles.

\begin{lemma}\label{lem:Xi}
Let $G$, $v$, $u$, and $w$ be as described above, and let $X_1=X_1^{uvw}$, $X_2=X_2^{uvw}$, and $X_3=X_3^{uvw}$. Assume that there exists an induced 6-cycle in $G$ containing the path $uvw$. Then $G_1$, the bipartite subgraph of $G$ induced between $X_1$ and $X_2$, is acyclic and hence the number of induced 6-cycles in $G$ containing the path $uvw$ is at most $|X_3|(|X_1|+|X_2|-1)$. Similarly $G_2$, the bipartite subgraph of $G$ induced between $X_2$ and $X_3$, is acyclic and the number of induced 6-cycles in $G$ containing the path $uvw$ is at most $|X_1|(|X_2|+|X_3|-1)$.
\end{lemma}
\begin{proof}
Suppose for a contradiction that $G_1$ contains the cycle $x_1y_1\dots x_ky_k$ for some $k\geq 2$, where $x_i\in X_1$ and $y_i\in X_2$ for all $i$. Let $z_1\in X_3$ be a common neighbour of $w$ and~$y_1$, and let $z_2\in X_3$ be a common neighbour of $w$ and $y_2$. If $z_1\neq z_2$, then $\{x_1,x_2,w\}$ and $\{y_1,y_2,u\}$ form the partite sets of a subdivision of $K_{3,3}$ in $G$, which contradicts Kuratowski's theorem~\cite{kuratowski} (see also~\cite{diestel}). Otherwise $\{x_1,x_2,z_1\}$ and $\{y_1,y_2,u\}$ form such partite sets, again giving a contradiction. Hence $G_1$ is acyclic, as is $G_2$, similarly.

Every induced 6-cycle in $G$ containing the path $uvw$ contains an edge of $G_1$ and a vertex in $X_3$, and this edge and vertex pair determines the 6-cycle. Since $G_1$ is acyclic and $X_1$ and $X_2$ are non-empty, there are at most $|X_1|+|X_2|-1$ edges in $G_1$. Hence there are at most $|X_3|(|X_1|+|X_2|-1)$ induced 6-cycles in $G$ which contain the path $uvw$. Similarly, considering $G_2$ in place of $G_1$ and $X_1$ in place of $X_3$, $G$ contains at most $|X_1|(|X_2|+|X_3|-1)$ such 6-cycles. This completes the proof of Lemma~\ref{lem:Xi}.
\end{proof}

When we apply this result in Section~\ref{sec:vtx_in_few}, we will use the following basic facts.

\begin{fact}\label{fact_top}
\begin{enumerate}[label=(\alph*)]
\item\label{fact:1} If $a,b,c\geq 0$ and $ab\geq c^2$, then $a+b\geq 2c$ with equality if and only if $a=b=c$.
\item\label{fact:2} If $a,b,c,d\geq 0$, $ab\geq cd+1$, and $|c-d|\leq 1$, then $a+b>c+d$.
\end{enumerate}
\end{fact}

Let $G$ be a plane graph with a (not necessarily induced) $K_{2,7}$ subgraph. The drawing of this $K_{2,7}$ subgraph induced by the drawing of $G$ splits the plane into seven regions, six bounded and one unbounded. Following the authors of~\cite{ghosh2021maximum}, we say that this $K_{2,7}$ subgraph of $G$ is an \emph{empty $K_{2,7}$} in $G$ if the interiors of all six of the bounded regions formed by the $K_{2,7}$ subgraph contain no vertices of $G$. If $G$ contains such a copy of $K_{2,7}$ then we say that $G$ \emph{contains an empty $K_{2,7}$}.

Inspired by a definition in~\cite{savery2021planar}, we say that two vertices in a planar graph $G$ are \emph{principal neighbours} if they are adjacent and there is an induced 6-cycle in $G$ containing both of them. Adapting an observation by the authors of~\cite{ghosh2021maximum}, we note that given an empty $K_{2,7}$ in a plane graph $G$, one of the vertices in the part of size 7 (the `central' one, in the natural sense) has only two principal neighbours in $G$, namely the vertices in the part of the $K_{2,7}$ of size 2.

Indeed, if $u$ and $w$ are the vertices in the part of size 2, and $v_1,\dots,v_7$ are the vertices in the part of size 7 (labeled in a natural order so that $v_4$ is the central one), then $v_4$ takes its neighbours from among $u$, $w$, $v_3$, and $v_5$. Clearly no induced 6-cycle contains the path $uv_4v_3$ since $u$ and $v_3$ are neighbours, and the same applies for the paths $uv_4v_5$, $wv_4v_3$, and $wv_4v_5$. If there is an induced 6-cycle containing the path $v_3v_4v_5$ then there is a path of length 4 from $v_3$ to $v_5$ in $G$ which avoids $v_4$, $u$, and $w$. This is impossible due to the other common neighbours of $u$ and $w$ and thus the only principal neighbours of $v_4$ are $u$ and $w$.

The next lemma allows us to find empty $K_{2,7}$'s in plane graphs under certain conditions. The argument is based on an argument in~\cite{ghosh2021maximum}.

\begin{lemma}\label{lem:empty}
Let $0<c<1$ be a real constant and let $n$ be large relative to $c$. Let $G$ be an $n$-vertex plane graph and suppose that every vertex of $G$ is in at least $n^2/10$ induced 6-cycles. Suppose that $u$ and $w$ are distinct vertices of $G$ with a common neighbourhood of size at least $cn$. Then $G$ contains an empty $K_{2,7}$ whose part of size 2 is $\{u,w\}$ and hence there is a vertex in $G$ whose only principal neighbours are $u$ and $w$.
\end{lemma}
\begin{proof}
Let $t=|N(u)\cap N(w)|$, so $t\geq cn$. Label the vertices in $N(u)\cap N(w)$ as $v_1,\dots,v_t$ in a natural order, that is, such that the drawing of the complete bipartite graph with parts $\{u,w\}$ and $\{v_1,\dots,v_t\}$ induced by the drawing of $G$ splits the plane into $t$ regions $R_1,\dots,R_t$ where the boundary of $R_i$ is formed of the cycle $v_iuv_{i+1}w$ for all $i$ (with addition in the subscript taken modulo $t$), and where $R_t$ is unbounded.

Let $3\leq i\leq t-3$ and consider the region $R_i$. Suppose for a contradiction that there is at least one vertex of $G$ in the interior of $R_i$, but not more than $n^{1/3}$. Note that none of these vertices are in the common neighbourhood of $u$ and $w$. Let $z$ be a vertex in the interior of $R_i$. Then $z$ is in at least $n^2/10$ induced 6-cycles in $G$ by assumption.

We will obtain a contradiction by deriving an upper bound of less than $n^2/10$ on the number of induced 6-cycles containing $z$. First, the number of induced 6-cycles containing $z$ which only use vertices in the interior of $R_i$ or on the boundary of $R_i$ is at most $(n^{1/3}+3)^5\leq n^{7/4}$ since there are at most $n^{1/3}+3$ options for each of the vertices in the 6-cycle besides $z$. It remains to bound the number of induced 6-cycles containing~$z$ that also contain a vertex in the exterior of $R_i$.

Any such cycle contains at least two vertices on the boundary of $R_i$. Clearly it cannot contain all four vertices on the boundary since these form a 4-cycle. Moreover, it cannot contain three vertices on the boundary. Indeed, any such three vertices form a path in $G$, so since the 6-cycle is induced it must contain this path. Of the remaining three vertices in the 6-cycle, one is $z$ and another, $x$, is in the exterior of $R_i$. Since $z$ and $x$ are not adjacent, the final vertex in the 6-cycle is one of their common neighbours. This common neighbour must be the fourth vertex on the boundary of $R_i$, so the 6-cycle is not induced.

Hence any induced 6-cycle containing $z$ and a vertex in the exterior of $R_i$ contains exactly two vertices on the boundary of $R_i$, and these must be non-adjacent since they must appear non-consecutively in the 6-cycle. However, they cannot be $v_i$ and $v_{i+1}$ since there is no path of length at most 4 from $v_i$ to $v_{i+1}$ contained in the exterior of~$R_i$ because of the $t-2$ other common neighbours of $u$ and $w$. Therefore, any induced 6-cycle containing $z$ and a vertex in the exterior of $R_i$ contains $u$ and $w$.

Since $z$ is not a common neighbour of $u$ and $w$, of the three other vertices in the cycle exactly one or two must be in the interior of $R_i$. The number of such cycles in which exactly two of the other vertices are in the interior of $R_i$ is at most $n^{5/3}$ since there are at most $n^{1/3}$ options for each of the two other vertices in the interior of $R_i$, and at most~$n$ for the vertex in the exterior.

The number of such cycles with exactly one other vertex in the interior of $R_i$ is at most $3n^{4/3}$, since there are at most $n^{1/3}$ options for the other vertex in the interior of $R_i$, and there are at most $3n-6$ edges in the (planar) bipartite graph induced between $N(u)$ and $N(w)$, so the number of options for the edge contained in the exterior of $R_i$ is less than $3n$.

Therefore $z$ is contained at most $n^{7/4}+n^{5/3}+3n^{4/3}<n^2/10$ induced 6-cycles, which contradicts our assumption. So for $3\leq i\leq t-3$, if $R_i$ has a vertex in its interior, then there are least $n^{1/3}$ such vertices. Let $l=|\{3\leq i\leq t-3:R_i \text{ has a vertex in its interior}\}|$. Then $ln^{1/3}\leq n$, so $l\leq n^{2/3}$. Hence since $t\geq cn$, there exist six consecutive values of $i$ in $3,\dots,t-3$ such that $R_i$ has no vertices in its interior for each of these values of $i$. Thus,~$G$ contains an empty $K_{2,7}$ in which the part of size 2 is $\{u,w\}$, which proves the first part of the conclusion of Lemma~\ref{lem:empty}; the second part now follows from the discussion above.
\end{proof}

In Lemma~\ref{lem:good_cycle} and its subsequent corollary, we use Lemma~\ref{lem:empty} to show that if a large planar graph has no vertices in `few' induced 6-cycles, then it contains an induced 6-cycle reminiscent of those in the graphs in $\mathcal{F}_n$.

\begin{lemma}\label{lem:good_cycle}
Let $n$ be large and let $G$ be an $n$-vertex planar graph. Suppose that every vertex of $G$ is in at least $n^2/10$ induced 6-cycles. Suppose $v_1$ is a vertex of $G$ with at most five principal neighbours. Then $G$ contains an induced 6-cycle $u_1v_1u_2v_2u_3v_3$ such that the only principal neighbours of $v_2$ are $u_2$ and $u_3$, and the only principal neighbours of $v_3$ are~$u_1$ and $u_3$.
\end{lemma}

\begin{proof}
Let $u_1$ and $u_2$ be neighbours of $v_1$ such that there are at least $n^2/100$ induced 6-cycles in $G$ containing the path $u_1v_1u_2$. Let $X_1=X_1^{u_1v_1u_2}$, $X_2=X_2^{u_1v_1u_2}$, and $X_3=X_3^{u_1v_1u_2}$. Let $G_1$ and $G_2$ be the bipartite subgraphs of $G$ induced between partite classes~$X_1$ and $X_2$, and $X_2$ and $X_3$ respectively. By Lemma~\ref{lem:Xi}, $G_1$ and $G_2$ are acyclic. 

Every induced 6-cycle in $G$ containing $u_1v_1u_2$ contains an edge of $G_2$, and every edge in $G_2$ is in at most $n$ such cycles, since these cycles are determined by the vertex in $X_1$ that they contain. There are at most $n$ edges in $G_2$ since it is acyclic. Let $l$ be the number of edges in $G_2$ which are in at least $n/1000$ induced 6-cycles in $G$ containing the path $u_1v_1u_2$. Then 
\[
\frac{n^2}{100}\leq nl+\frac{(n-l)n}{1000},
\]
so $l\geq n/111$.

If $x_2\in X_2$ and $x_3\in X_3$ are neighbours and there are at least $n/1000$ induced 6-cycles in $G$ containing the path $u_1v_1u_2x_3x_2$, then $x_2$ has at least $n/1000$ neighbours in $X_1$. Hence each of the $l$ edges above contains a vertex in $X_2$ which has degree at least $n/1000$ in $G_1$. Since $G_1$ is acyclic it contains at most $n$ edges, so at most 1000 vertices in $X_2$ can have degree at least $n/1000$ in $G_1$. Thus there exists $u_3\in X_2$ which is an endpoint of at least $n/111000$ of the $l$ edges, and this $u_3$ satisfies $|N(u_1)\cap N(u_3)|\geq n/1000$ and $|N(u_2)\cap N(u_3)|\geq n/111000$.

By Lemma~\ref{lem:empty}, there exist vertices $v_3\in N(u_1)\cap N(u_3)$ and $v_2\in N(u_2)\cap N(u_3)$ whose only principal neighbours are $u_1$ and $u_3$, and $u_2$ and $u_3$ respectively. Moreover, the same lemma implies that we can pick $v_2$ so that its only neighbours are $u_2$, $u_3$, and some vertices in the neighbourhood of $u_2$. Similarly, we can pick $v_3$ such that its only neighbours are $u_1$, $u_3$, and some vertices in the neighbourhood of $u_1$. Since $u_1$ and $u_2$ are not neighbours, we have $v_3\neq u_2,u_3$, and $v_3\not\in N(u_2)$. Hence $v_2$ and $v_3$ are not neighbours. The cycle $u_1v_1u_2v_2u_3v_3$ is therefore an induced 6-cycle in $G$, which completes the proof of Lemma~\ref{lem:good_cycle}.
\end{proof}

\begin{corollary}\label{cor:good_cycle}
Let $n$ be large and let $G$ be an $n$-vertex planar graph. Suppose that every vertex of $G$ is in at least $n^2/10$ induced 6-cycles. Then $G$ contains an induced 6-cycle $u_1v_1u_2v_2u_3v_3$ such that for each $i\in \{1,2,3\}$ the only principal neighbours of $v_i$ are $u_i$ and $u_{i+1}$, where addition in the subscript is taken modulo 3.
\end{corollary}

\begin{proof}
Every planar graph contains a vertex of degree at most 5. Applying Lemma~\ref{lem:good_cycle} to such a vertex of $G$ shows that there exists a vertex $v_1$ in $G$ with only two principal neighbours. Applying Lemma~\ref{lem:good_cycle} to $v_1$ completes the proof of Corollary~\ref{cor:good_cycle}.
\end{proof}

\section{Large planar graphs contain a vertex in few induced 6-cycles}\label{sec:vtx_in_few}

In this section we will show that every large planar graph contains a vertex in `few' induced 6-cycles, and that if the number of vertices in the graph is a multiple of 3 and no vertex is in `very few' induced 6-cycles, then the graph is in the family $\mathcal{F}_n$ for the appropriate $n$. The next two lemmas deal with the case where $n$ is a multiple of 3.

\begin{lemma}\label{lem:large_cmn_nbhd}
Let $n\equiv 0\pmod{3}$ be large, and let $G$ be an $n$-vertex planar graph. Suppose that every vertex of $G$ is contained in at least $\left(\frac{n}{3}-1\right)^2$ induced 6-cycles. If $G$ contains a (not necessarily induced) copy of $K_{2,n/3-1}$, then $G\in \mathcal{F}_n$.
\end{lemma}
\begin{proof}
Suppose that $G$ contains a copy of $K_{2,n/3-1}$. Fix such a copy and let $u_1$ and $u_2$ be the vertices in the part of size 2. By Lemma~\ref{lem:empty}, $G$ contains a vertex $v_1$ whose only principal neighbours are $u_1$ and $u_2$, so by Lemma~\ref{lem:good_cycle}, $G$ contains an induced 6-cycle $u_1v_1u_2v_2u_3v_3$ such that for each $i\in \{1,2,3\}$ the only principal neighbours of $v_i$ are $u_i$ and $u_{i+1}$, with subscript addition taken modulo 3.

Let $A=N(u_1)\cap N(u_2)$, so $|A|\geq \frac{n}{3}-1$. Now let $X=X^{u_1v_1u_2}$, $X_1=X_1^{u_1v_1u_2}$, $X_2=X_2^{u_1v_1u_2}$, and $X_3=X_3^{u_1v_1u_2}$. Note that $|X|\leq \frac{2n}{3}-1$ since $A$, $X$, and $\{u_1,u_2\}$ are  pairwise disjoint. By Lemma~\ref{lem:Xi} we have 
\begin{equation}\label{eq:large_cmn_nbhd_0}
    |X_3|(|X_1|+|X_2|-1)\geq \left(\frac{n}{3}-1\right)^2
\end{equation}
and
\begin{equation}\label{eq:large_cmn_nbhd_1}
    |X_1|(|X_2|+|X_3|-1)\geq \left(\frac{n}{3}-1\right)^2.
\end{equation}
Applying Fact~\hyperref[fact:1]{\ref*{fact_top}\ref*{fact:1}} gives $|X|\geq \frac{2n}{3}-1$, so in fact $|X|= \frac{2n}{3}-1$. This implies that $|A|=\frac{n}{3}-1$.  Applying the equality statement of Fact~\hyperref[fact:1]{\ref*{fact_top}\ref*{fact:1}} to inequality~\eqref{eq:large_cmn_nbhd_0} gives $|X_3|=\frac{n}{3}-1$, and applying it to inequality~\eqref{eq:large_cmn_nbhd_1} gives $|X_1|=\frac{n}{3}-1$. In turn, this implies $|X_2|=1$, so $X_2=\{u_3\}$. Let $B=X_3$ and $C=X_1$. Then $\{u_1,u_2,u_3\}$, $A$, $B$, and $C$ form a partition of the vertex set of $G$, and for each $a\in A$, $b\in B$, and $c\in C$ the edges $au_1$, $au_2$, $bu_2$, $bu_3$, $cu_3$ and $cu_1$ are present in $G$.

Since $u_1v_1u_2v_2u_3v_3$ is an induced 6-cycle in $G$, there are no edges between any of $u_1$, $u_2$, and $u_3$. Inspecting the argument of Lemma~\ref{lem:Xi} and noting that inequalities~\eqref{eq:large_cmn_nbhd_0} and~\eqref{eq:large_cmn_nbhd_1} are in fact equalities, we deduce that there are no edges between $B$ and $C$. By the definition of $X$ there are no edges between $u_1$ and $B$ or between $u_2$ and $C$. Repeating the arguments above with each of $v_2$ and $v_3$ in place of $v_1$ shows that there are no edges between $u_3$ and $A$, between $A$ and $B$, or between $A$ and $C$.

The only other possible edges of $G$ are edges within each of the vertex classes $A$, $B$, and $C$. By the planarity of $G$, for each of these vertex classes any edges within that vertex class clearly form a subset of the edge set of some path through the vertex set. Hence $G\in \mathcal{F}_n$, and Lemma~\ref{lem:large_cmn_nbhd} is proved.
\end{proof}

\begin{lemma}\label{lem:n_equiv_0}
Let $n\equiv 0\pmod{3}$ be large. Let $G$ be an $n$-vertex planar graph and suppose that every vertex of $G$ is in at least $\left(\frac{n}{3}-1\right)^2$ induced 6-cycles. Then $G\in \mathcal{F}_n$.
\end{lemma}
\begin{proof}
By Corollary~\ref{cor:good_cycle}, $G$ contains an induced 6-cycle $u_1v_1u_2v_2u_3v_3$ such that for each $i\in \{1,2,3\}$ the only principal neighbours of $v_i$ are $u_i$ and $u_{i+1}$, where addition in the subscript is taken modulo 3. Let $X=X^{u_1v_1u_2}$, $X_1=X_1^{u_1v_1u_2}$, $X_2=X_2^{u_1v_1u_2}$, $X_3=X_3^{u_1v_1u_2}$, $Y=X^{u_2v_2u_3}$, and $Z=X^{u_3v_3u_1}$.

By Lemma~\ref{lem:Xi} and Fact~\hyperref[fact:1]{\ref*{fact_top}\ref*{fact:1}} we have $|X|\geq \frac{2n}{3}-1$, and similarly $|Y|,|Z|\geq \frac{2n}{3}-1$. Suppose $|X|,|Y|,|Z|\geq \frac{2n}{3}+2$. Then $|X|+|Y|+|Z|\geq 2n+6$. Each of $u_1$, $u_2$, and $u_3$ appears in exactly one of $X$, $Y$, and $Z$, so $|X\setminus\{u_1,u_2,u_3\}|+|Y\setminus\{u_1,u_2,u_3\}|+|Z\setminus\{u_1,u_2,u_3\}|\geq 2n+3=2(n-3)+9$. So there are at least nine vertices in $G$ which appear in all three of $X$, $Y$, and $Z$. The following claim asserts that this is impossible.

\begin{claim}\label{claim:8_in_int}
There are at most eight vertices in $X\cap Y\cap Z$.
\end{claim}
\begin{proof}
If $z\in X\cap Y\cap Z$, then since $z$ is in an induced 6-cycle in $G$ containing $v_1$, $z$ is not in the common neighbourhood of $u_1$ and $u_2$. Similarly $z$ is not in the common neighbourhood of $u_1$ and $u_3$ or $u_2$ and $u_3$. So every vertex in $X\cap Y\cap Z$ is adjacent to at most one of $u_1$, $u_2$, and $u_3$.

Let $\Gamma$ be the induced 6-cycle $u_1v_1u_2v_2u_3v_3$. Fix a drawing of $G$ and note that in this drawing $\Gamma$ splits the plane into two regions. We will show that in the interiors of each of these regions there can be at most four vertices in $X\cap Y\cap Z$. Since none of the vertices in $\Gamma$ are in $X\cap Y\cap Z$ this will be sufficient to prove the claim.

First suppose that $z,z'\in X\cap Y\cap Z\cap N(u_1)$ are distinct vertices in one of the regions formed by $\Gamma$. Since $z$ is in an induced 6-cycle containing $v_2$, it has neighbours $a$ and $b$ such that $zau_2v_2u_3b$ is an induced 6-cycle. Note that $a$ and $b$ are distinct from $u_1$ since they are adjacent to $u_2$ and $u_3$ respectively, and they are distinct from each of $v_1$ and $v_3$ because they are in an induced 6-cycle which does not contain $u_1$.

Let $H_0$ be the subgraph of $G$ consisting of the cycle $\Gamma$ and the edges $u_1z$, $za$, $au_2$, $zb$, and $bu_3$, as shown in Figure~\hyperref[fig:8_in_int1]{\ref*{fig:8_in_int}(\subref*{fig:8_in_int1})}. The drawing of $H_0$ induced by the drawing of $G$ splits the plane into four regions, whose boundaries are formed of the cycles $\Gamma$, $u_1zau_2v_1$, $u_1zbu_3v_3$, and $zau_2v_2u_3b$. Since $a\in N(u_2)$ and $b\in N(u_3)$, $z'$ is distinct from both of these vertices. Hence we may assume that $z'$ is in the region with boundary $u_1zau_2v_1$. Since $z'$ is in an induced 6-cycle containing $v_2$, it has a common neighbour with $u_3$, and this common neighbour must be in the cycle $u_1zau_2v_1$. However, none of these vertices are adjacent to $u_3$, so we have a contradiction. Hence there are at most two vertices in $X\cap Y\cap Z\cap N(u_1)$, with at most one in each of the regions formed by $\Gamma$. Similarly there are at most four vertices in $X\cap Y\cap Z$ which are adjacent to $u_2$ or $u_3$.

\begin{figure}
    \centering
    \begin{subfigure}[h]{0.3\textwidth}
        \centering
        \includegraphics[width=\textwidth]{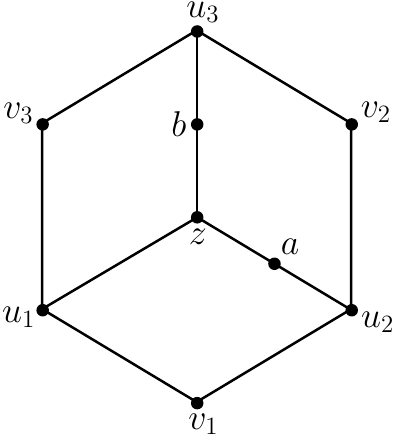}
        \caption{Graph $H_0$}
        \label{fig:8_in_int1}
     \end{subfigure}
     \hfill
     \begin{subfigure}[h]{0.3\textwidth}
        \centering
        \includegraphics[width=\textwidth]{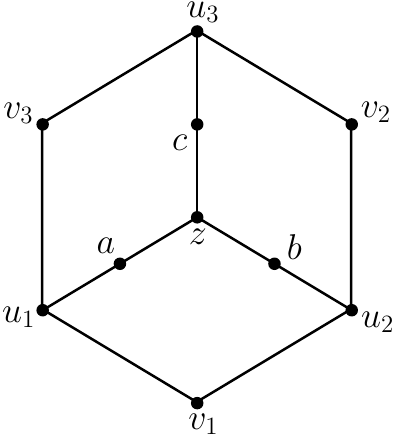}
        \caption{Graph $H_1$}
        \label{fig:8_in_int2}
     \end{subfigure}
    \hfill
     \begin{subfigure}[h]{0.3\textwidth}
        \centering
        \includegraphics[width=\textwidth]{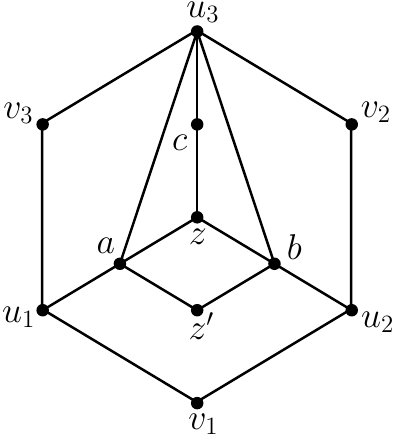}
        \caption{Graph $H_2$}
        \label{fig:8_in_int3}
     \end{subfigure}
        \caption{Some graphs used in the proof of Claim~\ref{claim:8_in_int}}
        \label{fig:8_in_int}
\end{figure}

Now suppose that $z,z'\in X\cap Y\cap Z$ are distinct vertices in one of the regions formed by $\Gamma$ which have no neighbours among $u_1$, $u_2$, and $u_3$. Note that $z$ and $z'$ are not adjacent to any of $v_1$, $v_2$, or $v_3$ since these vertices only have two principal neighbours, but they each share an induced 6-cycle with $z$ and $z'$. Since $z$ is in an induced 6-cycle containing~$v_1$, it has neighbours $x_1$ and $x_2$ such that $v_1u_1x_1zx_2u_2$ is an induced 6-cycle. 

Suppose first that $x_1\in N(u_3)$. Then since $z$ is in an induced 6-cycle containing~$v_3$ it has neighbours $x_3$ and $x_4$ such that $v_3u_3x_3zx_4u_1$ is an induced 6-cycle. Both $x_3$ and~$x_4$ are distinct from $x_1$, and at most one of them can be equal to $x_2$, so there exist distinct $a,b,c\in N(z)$ such that $a\in N(u_1)$, $b\in N(u_2)$, and $c\in N(u_3)$. If instead $x_1\not\in N(u_3)$, then since $z$ is in an induced 6-cycle containing $v_2$, there exist $x_5$ and $x_6$ such that $v_2u_2x_5zx_6u_3$ is an induced 6-cycle. Neither $x_5$ nor $x_6$ can be equal to $x_1$ since $x_1\not\in N(u_2)\cup N(u_3)$, so again there exist distinct $a$, $b$, and $c$ as above.

Fix such $a$, $b$, and $c$, and let $H_1$ be the subgraph of $G$ consisting of the cycle $\Gamma$ and the edges $za$, $zb$, $zc$, $au_1$, $bu_2$, and $cu_3$, as shown in Figure~\hyperref[fig:8_in_int2]{\ref*{fig:8_in_int}(\subref*{fig:8_in_int2})}. The drawing of $H_1$ induced by the drawing of $G$ splits the plane into four regions with boundaries formed of the cycles $\Gamma$, $zau_1v_1u_2b$, $zbu_2v_2u_3c$, and $zcu_3v_3u_1a$.

Since $z'$ is not in $H_1$, we may assume that $z'$ is in the region with boundary $zau_1v_1u_2b$. As $z'$ is in an induced 6-cycle containing $v_2$, $z'$ must have a common neighbour with $u_3$. This common neighbour must be in the cycle $zau_1v_1u_2b$ so we may assume that it is $a$. There is an induced 6-cycle in $G$ containing $z'$ and $v_3$, so $z'$ has a neighbour which is adjacent to $u_3$ but not to $u_1$. This neighbour must be $b$. Let $H_2$ be the subgraph of $G$ obtained by adding the edges $z'a$, $au_3$, $z'b$, and $bu_3$ to $H_1$, as shown in Figure~\hyperref[fig:8_in_int3]{\ref*{fig:8_in_int}(\subref*{fig:8_in_int3})}.

Since $z$ is in an induced 6-cycle containing $v_3$, it has a common neighbour $a'$ with $u_1$ which is not a neighbour of $u_3$. The cycle $az'bu_3$ separates $z$ and $u_1$ in $G$, so $a'$ must be a vertex in this cycle. However $a$ and $b$ are both adjacent to $u_3$ and neither $u_3$ nor $z'$ are adjacent to $u_1$, so no such $a'$ exists. This gives the desired contradiction, and hence there are at most two vertices in $X\cap Y\cap Z$ which are adjacent to none of $u_1$, $u_2$, and $u_3$. This completes the proof of Claim~\ref{claim:8_in_int}.
\end{proof}

Hence we may assume that $|X|\in \{\frac{2n}{3}-1, \frac{2n}{3}, \frac{2n}{3}+1\}$. We will consider each of these cases in turn. Note that if a vertex in $X_2$ has $\frac{n}{3}-1$ neighbours in $X_3$, then these are all common neighbours with $u_2$, so by Lemma~\ref{lem:large_cmn_nbhd}, $G\in \mathcal{F}_n$. Recall also that by the definitions of $X_2$ and $X_3$, every vertex in $X_2$ has a neighbour in $X_3$ and vice versa. Suppose first that $|X| = \frac{2n}{3}-1$. Then by Lemma~\ref{lem:Xi} and Fact~\hyperref[fact:1]{\ref*{fact_top}\ref*{fact:1}}, we have $|X_1|=|X_3|=\frac{n}{3}-1$. Hence $|X_2|=1$, so the lone vertex in $X_2$ has at least $\frac{n}{3}-1$ neighbours in $X_3$ which implies that $G\in\mathcal{F}_n$.

For the remaining two possible values of $|X|$ we will make use of the following claim. As in Lemmas~\ref{lem:Xi} and~\ref{lem:good_cycle}, let $G_1$ be the bipartite subgraph of $G$ induced between $X_1$ and~$X_2$, and let $G_2$ be that induced between $X_2$ and $X_3$.

\begin{claim}\label{claim:method1}
Let $s\leq 3$ be a positive integer. Suppose that there exists a set $T$ of at most four vertices in $X_2$ such that after excluding any one of them, those that remain have total degree at least $s$ in $G_2$. Then $s\leq |X|-\frac{2n}{3}+1$, and if $s=|X|-\frac{2n}{3}+1$, then $G_2$ is connected.
\end{claim}
\begin{proof}
Both $G_1$ and $G_2$ are acyclic by Lemma~\ref{lem:Xi}, and by the proof of that lemma an induced 6-cycle in $G$ containing $v_1$ is determined by the vertex in $X_1$ and the edge of $G_2$ that it contains. If $v\in X_1$ has at most one neighbour in $T$, then there are at least $s$ edges in $G_2$ which are not in an induced 6-cycle in $G$ containing $v$ and $v_1$. Since $G_1$ is acyclic, any pair of distinct vertices in $X_2$ has at most one common neighbour in $X_1$, so all but at most six vertices in $X_1$, have at most one neighbour in $T$.

Thus the number of induced 6-cycles in $G$ containing $v_1$ is at most $|X_1|(|X_2|+|X_3|-1)-s(|X_1|-6)=|X_1|(|X_2|+|X_3|-(s+1))+6s$. So $|X_1|(|X_2|+|X_3|-(s+1))\geq \left(\frac{n}{3}-1\right)^2-6s\geq \left(\frac{n-4}{3}\right)^2$ since $n$ is large. Noting that $G_2$ has at least $s+1$ edges so $|X_2|+|X_3|-(s+1)\geq 1$, we can apply Fact~\hyperref[fact:1]{\ref*{fact_top}\ref*{fact:1}} to obtain $|X|-(s+1)\geq \frac{2n-8}{3}$. Since the left-hand side is an integer, this implies that $|X|-(s+1)\geq \frac{2n}{3}-2$, or equivalently $s\leq |X|-\frac{2n}{3}+1$. If moreover $G_2$ is not connected, then the number of induced 6-cycles containing $v_1$ is at most $|X_1|(|X_2|+|X_3|-2)-s(|X_1|-6)$ which gives $s\leq |X|-\frac{2n}{3}$ by a similar argument. This completes the proof of Claim~\ref{claim:method1}.
\end{proof}

Suppose that $|X|=\frac{2n}{3}$. Then $|X|-\frac{2n}{3}+1=1$ so by Claim~\ref{claim:method1}, $|X_2|\leq 2$. If $|X_2|=1$, then we may assume that $|X_3|\geq \frac{n}{3}$, so the lone vertex in $X_2$ has at least $\frac{n}{3}$ neighbours in $X_3$ which implies $G\in \mathcal{F}_n$. If $|X_2|= 2$, then we may assume that $|X_3|\geq \frac{n}{3}-1$. By Claim~\ref{claim:method1} one vertex in $X_2$ has only one neighbour in $X_3$ and $G_2$ is connected, so the other vertex in $X_2$ has at least $\frac{n}{3}-1$ neighbours in $X_3$ and $G\in \mathcal{F}_n$.

Now suppose that $|X|=\frac{2n}{3}+1$. Then $|X|-\frac{2n}{3}+1=2$, so by Claim~\ref{claim:method1}, $|X_2|\leq 3$. If $|X_2|=1$, then we may assume that $|X_3|\geq \frac{n}{3}$, and so the lone vertex in $X_2$ has at least $\frac{n}{3}$ neighbours in $X_3$ and $G\in \mathcal{F}_n$. If $|X_2|=2$, then we may assume that $|X_3|\geq \frac{n}{3}$. By the claim, one vertex in $X_2$ has degree at most 2 in $G_2$. Moreover, if one vertex in $X_2$ has degree exactly 2 in $G_2$, then $G_2$ is connected so the other has at least $\frac{n}{3}-1$ neighbours in $X_3$ and $G\in\mathcal{F}_n$. If instead one vertex in $X_2$ has exactly one neighbour in $X_3$, then again the other has at least $\frac{n}{3}-1$, and $G\in\mathcal{F}_n$. Finally if $|X_2|=3$, then we may assume that $|X_3|\geq \frac{n}{3}-1$. By the claim, some pair of vertices in $X_2$ each have only one neighbour in $X_3$ and $G_2$ is connected, so the remaining vertex has at least $\frac{n}{3}-1$ neighbours in $X_3$. Thus $G\in \mathcal{F}_n$. This completes the proof of Lemma~\ref{lem:n_equiv_0}.
\end{proof}

The next two lemmas concern the cases where $n\equiv 1\pmod{3}$ or $n\equiv 2 \pmod{3}$. Their proofs make heavy use of some of the ideas in the previous two proofs in this section.

\begin{lemma}\label{lem:n_equiv_1}
Let $n$ be large and assume $n\equiv 1\pmod{3}$. Then every $n$-vertex planar graph contains a vertex in at most $\left(\frac{n-4}{3}\right)^2$ induced 6-cycles.
\end{lemma}
\begin{proof}
Let $G$ be an $n$-vertex planar graph and suppose for a contradiction that every vertex in $G$ is in at least $\left(\frac{n-4}{3}\right)^2+1$ induced 6-cycles. Then by Corollary~\ref{cor:good_cycle}, $G$ contains a vertex with only two principal neighbours. Let $v_1$ be such a vertex and apply Lemma~\ref{lem:good_cycle} to obtain an induced 6-cycle $u_1v_1u_2v_2u_3v_3$ such that for each $i\in \{1,2,3\}$ the only principal neighbours of $v_i$ are $u_i$ and $u_{i+1}$, where addition in the subscript is taken modulo 3.

Let $X$, $X_1$, $X_2$, $X_3$, $Y$, $Z$, $G_1$, and $G_2$ be as in the proof of Lemma~\ref{lem:n_equiv_0}. By Lemma~\ref{lem:Xi} and Fact~\hyperref[fact:2]{\ref*{fact_top}\ref*{fact:2}} we have $|X|\geq \frac{2n-2}{3}$, and similarly $|Y|,|Z|\geq \frac{2n-2}{3}$. Claim~\ref{claim:8_in_int} still holds in this setting by the same proof, and applying this in the same way as in the proof of Lemma~\ref{lem:n_equiv_0}, we deduce that at least one of $X$, $Y$, and $Z$ has size less than $\frac{2n+7}{3}$. Hence we may assume that $|X|\in \{\frac{2n-2}{3}, \frac{2n+1}{3}, \frac{2n+4}{3}\}$. We will shortly consider each of these possible values of $|X|$ in turn and show that each leads to a contradiction.

First, we state a version of Claim~\ref{claim:method1} for this setting.

\begin{claim}\label{claim:method2}
    Let $s\leq 4$ be a positive integer. Suppose that there exists a set of at most five vertices in $X_2$ such that after excluding any one of them, those that remain have total degree at least $s$ in $G_2$. Then $s\leq |X|-\frac{2n-5}{3}$, and if $s=|X|-\frac{2n-5}{3}$, then $G_2$ is connected. 
\end{claim}

It is straightforward to modify the proof of Claim~\ref{claim:method1} to show Claim~\ref{claim:method2}. We will also need the following. 

\begin{claim}\label{claim:cmn_nbhd1}
Suppose that for every triple of distinct vertices $u$, $v$, and $w$ in $G$ such that $u$ and $w$ are the only principal neighbours of $v$ we have $|X^{uvw}|\geq \frac{2n+t}{3}$ where $t\in \{-2,1,4\}$. Then $G$ contains no copy of $K_{2,(n-t)/3-1}$.
\end{claim}
\begin{proof}
Suppose for a contradiction that $G$ contains a copy of $K_{2,(n-t)/3-1}$. Let $u$ and $w$ be the vertices in the part of this bipartite graph of size 2. By Lemma~\ref{lem:empty}, $G$ contains a vertex $v$ whose only principal neighbours are $u$ and $w$, and by assumption $|X^{uvw}|\geq \frac{2n+t}{3}$. However, $X^{uvw}$ is disjoint from the vertex set of the copy of $K_{2,(n-t)/3-1}$, which contradicts the fact that $G$ contains only $n$ vertices. This completes the proof of Claim~\ref{claim:cmn_nbhd1}.
\end{proof}

Thus $G$ contains no copy of $K_{2,(n-1)/3}$. We now consider each of the three possible values of $|X|$. Suppose first that $|X|=\frac{2n-2}{3}$, then $|X|-\frac{2n-5}{3}=1$ so by Claim~\ref{claim:method2}, $|X_2|\leq 2$. If $|X_2|=1$, then we may assume that $|X_3|\geq \frac{n-1}{3}$ so the single vertex in $X_2$ has at least $\frac{n-1}{3}$ common neighbours with $u_2$, which is a contradiction.

If $|X_2|=2$, then $G_2$ is connected and one vertex in $X_2$ has degree 1 in $G_2$. By symmetry, the same holds if we replace $G_2$ by $G_1$. If one vertex in $X_2$ has degree 1 in $G_1$ and the other has degree 1 in $G_2$, then the number of induced 6-cycles containing $v_1$ is clearly at most linear in $n$, so in fact one vertex in $X_2$ has degree 1 in $G_1$ and in $G_2$.

Thus one vertex in $X_2$ is in at most one induced 6-cycle containing $v_1$, and the other is in at most $|X_1||X_3|$. Hence
\begin{equation}\label{eq:special_case}
|X_1||X_3|+1\geq \left(\frac{n-4}{3}\right)^2+1,
\end{equation}
and $|X_1|+|X_3|=\frac{2n-8}{3}$ so by Fact~\hyperref[fact:1]{\ref*{fact_top}\ref*{fact:1}}, $|X_1|=|X_3|=\frac{n-4}{3}$. Moreover, since we have equality in~\eqref{eq:special_case}, there are no edges between $X_1$ and $X_3$. Note that $v_3\in X_1$ and $u_3\in X_2$. If $u_3$ is the vertex in $X_2$ of degree 1 in $G_1$, then~$v_3$ is its sole neighbour in $X_1$. But then~$v_3$ is in an induced 6-cycle in $G$ containing the other vertex in $X_2$, which contradicts the fact that $u_1$ and $u_3$ are its only principal neighbours. Hence $u_3$ is the vertex in $X_2$ of degree $\frac{n-4}{3}$ in $G_1$ and $G_2$. Let the other vertex in $X_2$ be $z$, and let $B=X_3$ and $C=X_1$. This is shown in Figure~\hyperref[fig:special_case1]{\ref*{fig:special_case}(\subref*{fig:special_case1})}.

\begin{figure}
    \centering
    \begin{subfigure}[h]{0.45\textwidth}
        \centering
        \includegraphics[width=\textwidth]{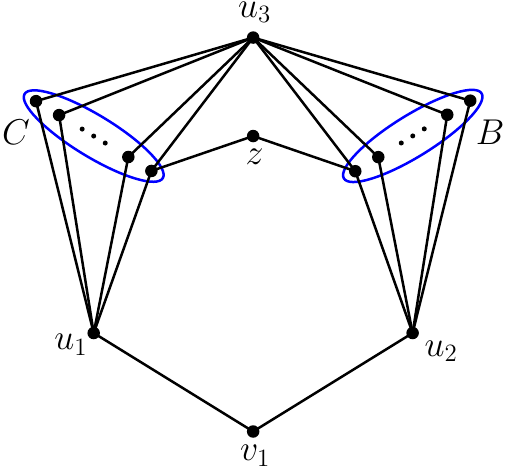}
        \caption{The induced 6-cycles containing $v_1$}
        \label{fig:special_case1}
     \end{subfigure}
     \hfill
     \begin{subfigure}[h]{0.45\textwidth}
        \centering
        \includegraphics[width=\textwidth]{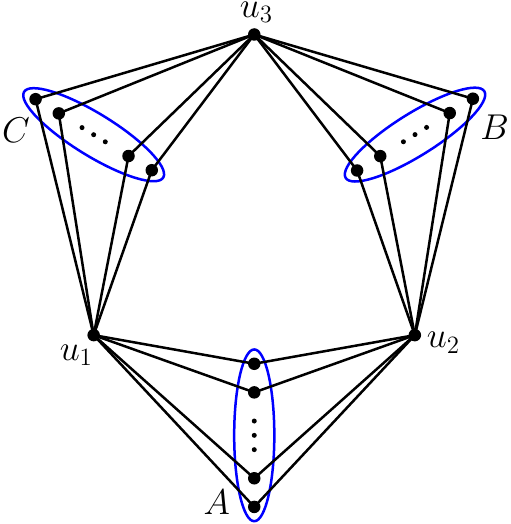}
        \caption{Graph $G'$}
        \label{fig:special_case2}
     \end{subfigure}
        \caption{Two graphs used in the proof of Lemma~\ref{lem:n_equiv_1}}
        \label{fig:special_case}
\end{figure}

Note that $B$ and $Y$ are disjoint, so $|Y|\leq \frac{2n-2}{3}$. Hence in fact $|Y| = \frac{2n-2}{3}$. Applying the arguments above with $Y$ in place of $X$ we deduce that $u_1$ and $u_2$ have at least $\frac{n-4}{3}$ common neighbours. Since $G$ has only $n$ vertices and $X\cap N(u_1)\cap N(u_2)=\emptyset$, $u_1$ and~$u_2$ must in fact have exactly $\frac{n-4}{3}$ common neighbours. Let $A$ be this set of common neighbours.

Let $G'$ be the subgraph of $G$ consisting of the edges $au_1$, $au_2$, $bu_2$, $bu_3$, $cu_3$, and $cu_1$ for each $a\in A$, $b\in B$, and $c\in C$, as shown in Figure~\hyperref[fig:special_case2]{\ref*{fig:special_case}(\subref*{fig:special_case2})}. Fix a drawing of $G$, and consider the induced drawing of $G'$. No face in the drawing of $G'$ has more than six vertices on its boundary, so $z$ has degree at most 6 in $G$. By Lemma~\ref{lem:Xi}, in order for $z$ to be in at least $\left(\frac{n-4}{3}\right)^2+1$ induced 6-cycles, it must have two distinct neighbours both of whose degrees in $G$ are linear in $n$. By considering the faces in the drawing of $G'$ again, we see that the only vertices in $G-\{z\}$ of degree greater than 6 are $u_1$, $u_2$, and $u_3$, but~$z$ is adjacent to neither $u_1$ nor $u_2$, so we have a contradiction.
Let $G'$ be the subgraph of $G$ consisting of the edges $au_1$, $au_2$, $bu_2$, $bu_3$, $cu_3$, and $cu_1$ for each $a\in A$, $b\in B$, and $c\in C$, as shown in Figure~\hyperref[fig:special_case2]{\ref*{fig:special_case}(\subref*{fig:special_case2})}. Fix a drawing of $G$, and consider the induced drawing of $G'$. No face in the drawing of $G'$ has more than six vertices on its boundary, so $z$ has degree at most 6 in $G$. By Lemma~\ref{lem:Xi}, in order for $z$ to be in at least $\left(\frac{n-4}{3}\right)^2+1$ induced 6-cycles, it must have two distinct neighbours both of whose degrees in $G$ are linear in $n$. By considering the faces in the drawing of $G'$ again, we see that the only vertices in $G-\{z\}$ of degree greater than 6 are $u_1$, $u_2$, and $u_3$, but $z$ is adjacent to neither $u_1$ nor $u_2$, so we have a contradiction.

Therefore $|X_2|\neq 2$ and hence $|X|\neq \frac{2n-2}{3}$. Applying Claim~\ref{claim:cmn_nbhd1} shows that $G$ therefore contains no copy of $K_{2,(n-4)/3}$. Next suppose that $|X|=\frac{2n+1}{3}$, so $|X|-\frac{2n-5}{3}=2$. Claim~\ref{claim:method2} implies that $|X_2|\leq 3$. If $|X_2|=1$, then we may assume that $|X_3|\geq \frac{n-1}{3}$ which gives rise to a copy of $K_{2,(n-1)/3}$ in $G$, and hence a contradiction. If $|X_2|=2$, then we may again assume that $|X_3|\geq \frac{n-1}{3}$. By Claim~\ref{claim:method2}, one vertex in $X_2$ has degree at most 2 in $G_2$, and if one has degree exactly 2, then $G_2$ is connected. In this case the other vertex has degree at least $\frac{n-4}{3}$ in $G_2$ which is a contradiction. If instead one vertex in $X_2$ has degree 1 in~$G_2$, then again the other has degree at least $\frac{n-4}{3}$, so again we arrive at a contradiction. If $|X_2|=3$, then we may assume that $|X_3|\geq \frac{n-4}{3}$. By Claim~\ref{claim:method2}, some two vertices in $X_2$ each have degree 1 in $G_2$ and $G_2$ is connected, so the remaining vertex in $X_2$ has at least $\frac{n-4}{3}$ neighbours in $X_3$ which is a contradiction.

Therefore $|X|\neq \frac{2n+1}{3}$, so $G$ contains no copy of $K_{2,(n-7)/3}$ by Claim~\ref{claim:cmn_nbhd1}. Suppose finally that $|X|=\frac{2n+4}{3}$, in which case $|X|-\frac{2n-5}{3}=3$ and Claim~\ref{claim:method2} implies that $|X_2|\leq 4$. If $|X_2|=1$, then we may assume that $|X_3|\geq \frac{n+2}{3}$ and hence $G$ has a copy of $K_{2,(n+2)/3}$ which is a contradiction. If $|X_2|=2$, then we may assume that $|X_3|\geq \frac{n-1}{3}$. By Claim~\ref{claim:method2}, one vertex in $X_2$ has degree at most 3 in $G_2$, and if one has degree exactly 3, then $G_2$ is connected. In this case, the other vertex in $X_2$ has at least $\frac{n-7}{3}$ neighbours in $X_3$, which is a contradiction. So one vertex in $X_2$ has at most two neighbours in $X_3$, and hence the other has at least $\frac{n-7}{3}$, which is another contradiction.

If $|X_2|=3$, then we may again assume that $|X_3|\geq \frac{n-1}{3}$. By Claim~\ref{claim:method2}, one vertex in $X_2$ has degree 1 in $G_2$, and another has degree at most 2. Moreover if this second vertex has degree exactly 2 in $G_2$, then $G_2$ is connected, so the remaining vertex in $X_2$ has at least $\frac{n-4}{3}$ neighbours in $X_3$, a contradiction. Otherwise some two vertices in $X_2$ each have only one neighbour in $X_3$, so the final vertex has at least $\frac{n-7}{3}$ such neighbours, a contradiction. If $|X_2|=4$, then we may assume that $|X_3|\geq \frac{n-4}{3}$. By Claim~\ref{claim:method2}, some three vertices in $X_2$ each have only one neighbour in $X_3$ and $G_2$ is connected, so the remaining vertex in $X_2$ has at least $\frac{n-4}{3}$ neighbours in $X_3$, which is a contradiction.

Thus $|X|\neq \frac{2n+4}{3}$, so $|X|\not\in\{\frac{2n-2}{3}, \frac{2n+1}{3}, \frac{2n+4}{3}\}$ which gives the required contradiction and completes the proof of Lemma~\ref{lem:n_equiv_1}.
\end{proof}

The proof of the next lemma is very similar to that of Lemma~\ref{lem:n_equiv_1}, but is more straightforward. We include it nonetheless, for completeness.

\begin{lemma}\label{lem:n_equiv_2}
Let $n$ be large and assume $n\equiv 2\pmod{3}$. Then every $n$-vertex planar graph contains a vertex in at most $\left(\frac{n-5}{3}\right)\left(\frac{n-2}{3}\right)$ induced 6-cycles.
\end{lemma}
\begin{proof}
Let $G$ be an $n$-vertex planar graph and suppose for a contradiction that every vertex in $G$ is in at at least $\left(\frac{n-5}{3}\right)\left(\frac{n-2}{3}\right)+1$ induced 6-cycles. Then by Corollary~\ref{cor:good_cycle}, $G$ contains a vertex with only two principal neighbours. Let $v_1$ be such a vertex and apply Lemma~\ref{lem:good_cycle} to obtain an induced 6-cycle $u_1v_1u_2v_2u_3v_3$ such that for each $i\in \{1,2,3\}$ the only principal neighbours of $v_i$ are $u_i$ and $u_{i+1}$ where addition in the subscript is taken modulo 3.

Let $X$, $X_1$, $X_2$, $X_3$, $Y$, $Z$, $G_1$, and $G_2$ be as in the proofs of Lemmas~\ref{lem:n_equiv_0} and~\ref{lem:n_equiv_1}. By Lemma~\ref{lem:Xi} and Fact~\hyperref[fact:2]{\ref*{fact_top}\ref*{fact:2}} we have $|X|\geq \frac{2n-1}{3}$, and similarly for $Y$ and $Z$. Again, Claim~\ref{claim:8_in_int} holds in this setting by the same argument, so applying this in the same way as before shows that at least one of $X$, $Y$, and $Z$ has size at most $\frac{2n+5}{3}$. Hence we may assume that $|X|\in \{\frac{2n-1}{3}, \frac{2n+2}{3}, \frac{2n+5}{3}\}$. As before, we will consider each of these cases in turn and will obtain a contradiction in each case.

The version of Claim~\ref{claim:method1} relevant to this setting (proved using Fact~\hyperref[fact:2]{\ref*{fact_top}\ref*{fact:2}} in the proof in place of Fact~\hyperref[fact:1]{\ref*{fact_top}\ref*{fact:1}}) is as follows.

\begin{claim}\label{claim:method3}
    Let $s\leq 4$ be a positive integer. Suppose that there exists a set of at most five vertices in $X_2$ such that after excluding any one of them, those that remain have total degree at least $s$ in $G_2$. Then $s\leq |X|-\frac{2n-4}{3}$, and if $s=|X|-\frac{2n-4}{3}$, then $G_2$ is connected
\end{claim}

The following version of Claim~\ref{claim:cmn_nbhd1} holds in this setting by exactly the same proof.

\begin{claim}\label{claim:cmn_nbhd2}
Suppose that for every triple of distinct vertices $u$, $v$, and $w$ in $G$ such that $u$ and $w$ are the only principal neighbours of $v$ we have $|X^{uvw}|\geq \frac{2n+t}{3}$ where $t\in \{-1,2,5\}$. Then $G$ contains no copy of $K_{2,(n-t)/3-1}$.
\end{claim}

By Claim~\ref{claim:cmn_nbhd2}, $G$ contains no copy of $K_{2,(n-2)/3}$. Suppose that $|X|= \frac{2n-1}{3}$. Then $|X|-\frac{2n-4}{3}=1$ so by the Claim~\ref{claim:method3}, $|X_2|\leq 2$. If $|X_2|=1$, then we may assume that $|X_3|\geq \frac{n-2}{3}$ which gives a copy of $K_{2,(n-2)/3}$ and hence a contradiction. If $|X_2|=2$, then we may again assume that $|X_3|\geq \frac{n-2}{3}$. By Claim~\ref{claim:method3}, one vertex in $X_2$ has degree 1 in $G_2$ and $G_2$ is connected, so the other vertex in $X_2$ has at least $\frac{n-2}{3}$ neighbours in $X_3$, which is a contradiction.

Hence $|X|\neq \frac{2n-1}{3}$, so by Claim~\ref{claim:cmn_nbhd2}, $G$ contains no copy of $K_{2,(n-5)/3}$. Suppose that $|X|=\frac{2n+2}{3}$, so $|X|-\frac{2n-4}{3}=2$ and hence by Claim~\ref{claim:method3}, $|X_2|\leq 3$. If $|X_2|=1$, then we may assume that $|X_3|\geq \frac{n+1}{3}$ which gives a copy of $K_{2,(n+1)/3}$ and hence a contradiction. If $|X_2|=2$, then we may assume that $|X_3|\geq \frac{n-2}{3}$. One vertex in $X_2$ has degree at most 2 in $G_2$, and if one has degree exactly 2 then $G_2$ is connected, so the other has degree at least $\frac{n-5}{3}$ which is a contradiction. Similarly if one vertex in $X_2$ has degree 1 in $G_2$, then the other has degree at least $\frac{n-5}{3}$. If $|X_2|= 3$, then again we may assume that $|X_3|\geq \frac{n-2}{3}$. By Claim~\ref{claim:method3} some two vertices in $X_2$ have degree 1 in $G_2$ and $G_2$ is connected, so the other vertex has degree at least $\frac{n-2}{3}$ in $G_2$, which is a contradiction.

Hence $|X|\neq \frac{2n+2}{3}$, so there is no copy of $K_{2,(n-8)/3}$ in $G$ by Claim~\ref{claim:cmn_nbhd2}. Suppose finally that $|X|=\frac{2n+5}{3}$, so $|X|-\frac{2n-4}{3}=3$ and hence $|X_2|\leq 4$. If $|X_2|=1$, then we may assume that $|X_3|\geq \frac{n+1}{3}$ which gives a copy of $K_{2,(n+1)/3}$ in $G$, and hence a contradiction. If $|X_2|=2$, then we may again assume that $|X_3|\geq \frac{n+1}{3}$. One vertex in $X_2$ has degree at most 3 in $G_2$, so the other has degree at least $\frac{n-8}{3}$ in $G_2$, which is a contradiction.

If $|X_2|= 3$, then we may assume that $|X_3|\geq \frac{n-2}{3}$. By Claim~\ref{claim:method3}, one vertex in $X_2$ has degree 1 in $G_2$ and another has degree at most 2. If this second vertex has degree exactly 2 in $G_2$, then $G_2$ is connected so the remaining vertex has degree at least $\frac{n-5}{3}$ in $G_2$, a contradiction. If instead two vertices in $X_2$ have degree 1 in $G_2$, then the remaining vertex has degree at least $\frac{n-8}{3}$, which is another contradiction. Finally, if $|X_2|=4$, then we may assume that $|X_3|\geq \frac{n-2}{3}$. Some three vertices in $X_2$ have degree 1 in $G_2$ and $G_2$ is connected, so the remaining vertex has degree at least $\frac{n-2}{3}$, which is a contradiction. Thus $|X|\neq \frac{2n+8}{3}$, so $|X|\not\in \{\frac{2n-1}{3}, \frac{2n+2}{3}, \frac{2n+5}{3}\}$ which gives the required contradiction and completes the proof of Lemma~\ref{lem:n_equiv_2}.
\end{proof}

\section{Proof of Theorem~\ref{thm:main}}\label{sec:main_proof}

We are now ready to prove Theorem~\ref{thm:main}.

\begin{proof}[Proof of Theorem~\ref{thm:main}]
First, clearly any graph in $\mathcal{F}_n$ is planar. For $n\geq 6$, let 
\[
h_0(n)=
\begin{cases}
\left(\frac{n}{3}-1\right)^3, & \text{if $n\equiv 0\pmod{3}$} \\
\left(\frac{n-4}{3}\right)^2\left(\frac{n-1}{3}\right), & \text{if $n\equiv 1\pmod{3}$} \\
\left(\frac{n-2}{3}\right)^2\left(\frac{n-5}{3}\right), & \text{if $n\equiv 2\pmod{3}$}.
\end{cases}
\]
Let $G\in \mathcal{F}_n$ and label its vertices as in Definition~\ref{def:Fn}. For each $a\in A$, $b\in B$, and $c\in C$, the cycle $u_1au_2bu_3c$ is an induced 6-cycle in $G$, so $G$ contains at least $h_0(n)$ induced 6-cycles.

Suppose there exists another induced 6-cycle, $\Gamma$, in $G$. We may assume that $\Gamma$ contains some $a\in A$. There is clearly no 6-cycle in $G$ consisting entirely of vertices in $A$, so we may also assume that $u_2$ is one of the neighbours of $a$ in $\Gamma$. If the other neighbour of $u_2$ in $\Gamma$ is some $b\in B$, then $\Gamma$ contains no further vertices in $B$ since they are all neighbours of $u_2$. Hence the other neighbour of $b$ in $\Gamma$ is $u_3$, and the other neighbour of $u_3$ in $\Gamma$ is some $c\in C$. The remaining vertex in $\Gamma$ is a common neighbour of $a$ and $c$, so it must be $u_1$. Hence $\Gamma$ is one of the induced 6-cycles described above, which is a contradiction. In the other case $\Gamma$ contains the path $au_2a'$ for some $a'\in A$. It contains no further vertices in $A$, so the other neighbour of $a'$ in $\Gamma$ is $u_1$, but this is also a neighbour of $a$, which contradicts $\Gamma$ being an induced 6-cycle. Therefore any graph in $\mathcal{F}_n$ contains exactly $h_0(n)$ induced 6-cycles.

For $n\geq 6$, define $\mathcal{G}_n$ to be the family of $n$-vertex planar graphs containing $f_I(n,C_6)$ induced 6-cycles. To prove the theorem it is sufficient to show that $\mathcal{G}_n\subseteq \mathcal{F}_n$ for large enough $n$. For $n\geq 6$ let 
\[
h_1(n)=
\begin{cases}
\left(\frac{n}{3}-1\right)^2, & \text{if $n\equiv 0\pmod{3}$} \\
\left(\frac{n-4}{3}\right)^2, & \text{if $n\equiv 1\pmod{3}$} \\
\left(\frac{n-2}{3}\right)\left(\frac{n-5}{3}\right), & \text{if $n\equiv 2\pmod{3}$}
\end{cases}
\]
so that $h_0(n)-h_0(n-1)=h_1(n)$ for all $n\geq 7$. By Lemmas~\ref{lem:n_equiv_1} and~\ref{lem:n_equiv_2}, when $n$ is large and not divisible by 3, every graph in $\mathcal{G}_n$ has a vertex in at most $h_1(n)$ induced 6-cycles. By taking such a graph and deleting a vertex in $h_1(n)$ induced 6-cycles, we obtain an $(n-1)$-vertex planar graph containing at least $f_I(n,C_6)-h_1(n)$ induced 6-cycles. Thus if $n$ is large and not divisible by 3, then $f_I(n,C_6)- f_I(n-1,C_6)\leq h_1(n)$.

By Lemma~\ref{lem:n_equiv_0}, when $n$ is large and divisible by 3, every graph in $\mathcal{G}_n$ is either in $\mathcal{F}_n$ or has a vertex in fewer than $h_1(n)$ induced 6-cycles. By the count above, every graph in $\mathcal{F}_n$ contains a vertex in exactly $h_1(n)$ induced 6-cycles, so by the same argument as in the case where $n$ is not a multiple of 3, if $n$ is large and divisible by 3, then $f_I(n,C_6)- f_I(n-1,C_6)\leq h_1(n)$.

If $f_I(n,C_6)- f_I(n-1,C_6) < h_1(n)$ for infinitely many values of $n$, then for large enough $n$ we have $h_0(n) > f_I(n,C_6)$, which is a contradiction, since every graph in $\mathcal{F}_n$ is planar and contains $h_0(n)$ induced 6-cycles. Therefore for large enough $n$, we have
\begin{equation}\label{eq:f_step}
    f_I(n,C_6) - f_I(n-1,C_6) = h_1(n).
\end{equation}
If $n$ is large enough that this equality holds and $G\in \mathcal{G}_n$, then every vertex in $G$ is in at least $h_1(n)$ induced 6-cycles, otherwise we could delete the vertex in the fewest such cycles to obtain an $(n-1)$-vertex planar graph containing more than $f_I(n-1,C_6)$ induced 6-cycles. Hence if $n$ is large and divisible by 3, and $G\in \mathcal{G}_n$, then $G\in \mathcal{F}_n$ by Lemma~\ref{lem:n_equiv_0}.

Now let $n$ be large with $n\equiv 1\pmod{3}$, and let $G\in \mathcal{G}_n$. By Lemma~\ref{lem:n_equiv_1} and equation~\eqref{eq:f_step} there exists a vertex in $G$ in exactly $h_1(n)$ induced 6-cycles. Deleting such a vertex, $z$, yields $G'\in \mathcal{G}_{n-1}$. Since $n-1\equiv 0 \pmod{3}$, we have $G'\in \mathcal{F}_{n-1}$. Label the vertices of $G'$ according to Definition~\ref{def:Fn}. If $z$ has more than two neighbours in any of $A$,~$B$, or $C$ in $G$, then this clearly gives rise to a subdivision of $K_{3,3}$, which contradicts the fact that $G$ is planar. Hence $z$ has degree at most 9 in $G$.

It follows that $z$ has two neighbours with which it is in at least $n^2/360$ induced 6-cycles. By Lemma~\ref{lem:Xi}, these two neighbours must each have degree at least $n/360$ in $G$, so without loss of generality $z$ is adjacent to $u_1$ and $u_2$ in $G$. Since $z$ is adjacent to at most two vertices in $A$, there exists $a\in A$ which has no path of length at most 4 to $z$ which avoids $u_1$ and $u_2$. Any induced 6-cycle in $G$ containing $a$ and $z$ must contain $u_1$ and $u_2$, but $au_1zu_2$ is a 4-cycle in $G$, so in fact $a$ and $z$ do not share an induced 6-cycle in $G$. Thus $a$ is in exactly $h_1(n-1)=h_1(n)$ induced 6-cycles in $G$.

Let $G''$ be the graph obtained by deleting $a$ from $G$, then as before we have $G''\in \mathcal{F}_{n-1}$. Since $z$ is adjacent to $u_1$ and $u_2$, it follows that $z$ is not adjacent to $u_3$, nor to any vertices in $B$ or $C$. So we obtain $G$ from $G'$ by adding a vertex whose neighbours are $u_1$, $u_2$, and some vertices in $A$. Since $G$ is planar, it is now clear that $G\in \mathcal{F}_n$.

Finally suppose $n\equiv 2\pmod{3}$ and let $G\in \mathcal{G}_n$. By Lemma~\ref{lem:n_equiv_2} and equation~\eqref{eq:f_step}, there exists a vertex, $z$, in $G$ in exactly $h_1(n)$ induced 6-cycles. Deleting this vertex yields $G'\in \mathcal{F}_{n-1}$. Label the vertices of $G'$ according to Definition~\ref{def:Fn}, with $A$ the larger of the three vertex classes. As before, $z$ is adjacent to at least two of $u_1$, $u_2$ and $u_3$.

If $z$ is adjacent to $u_1$ and $u_2$, then as in the case $n\equiv 1\pmod{3}$ there exists some $a\in A$ which is not in an induced 6-cycle containing $z$, and hence is in exactly $h_1(n-1)$ induced 6-cycles in $G$. But $h_1(n-1)< h_1(n)$, and every vertex of $G$ is in at least $h_1(n)$ induced 6-cycles by equation~\eqref{eq:f_step}, so we have a contradiction. Hence without loss of generality $z$ is adjacent to $u_2$ and $u_3$ in $G$.

Again, there is some $b\in B$ which is not in an induced 6-cycle in $G$ containing $z$, so $b$ is in exactly $h_1(n)$ induced 6-cycles in $G$. Deleting $b$ from $G$ yields another graph in $\mathcal{F}_{n-1}$, so it follows that $z$ has no neighbours in $A\cup C\cup \{u_1\}$, and thus we obtain $G$ from $G'$ by adding a vertex whose neighbours are $u_2$, $u_3$, and some vertices in $B$. By the planarity of $G$, it is clear that $G\in \mathcal{F}_n$. This completes the proof of Theorem~\ref{thm:main}.
\end{proof}

\section{Longer even cycles}\label{sec:longer_cycles}

In this section we discuss to what extent the ideas used above to prove Theorem~\ref{thm:main} might be of use in tackling Conjecture~\ref{conj:induced_even} for cycles of length greater than 6. In the case of odd cycles it seems that very little of the argument can be straightforwardly adapted, so we focus on the case of even cycles. Recall that we determined the value of $f_I(n,C_6)$ for large $n$ in roughly the following steps:
\begin{enumerate}
    \item For large $n$, we found some structure in those $n$-vertex planar graphs in which every vertex is in at least $n^2/10$ induced 6-cycles.
    \item For large $n$, we showed that if $G$ is an $n$-vertex planar graph, then it contains a vertex in at most $h_1(n)$ induced 6-cycles, where $h_1(n)$ is the minimum number of induced 6-cycles that a vertex in a graph in $\mathcal{F}_n$ is contained in.
    \item For $n$ large and divisible by 3 we showed that if $G$ is an $n$-vertex planar graph in which every vertex is in at least $h_1(n)$ induced 6-cycles, then $G\in\mathcal{F}_n$.
    \item We used the fact that in an $n$-vertex planar graph containing $f_I(n,C_6)$ induced 6-cycles every vertex is in at least $f_I(n,C_6)-f_I(n-1,C_6)$ induced 6-cycles, combined with step 2, to deduce that $f_I(n,C_6)=f_I(n-1,C_6)+h_1(n)$ for large enough $n$. We then used step 3 to determine the value of $f_I(n,C_6)$ exactly for large $n$. 
\end{enumerate}

In attempting to adapt this to $2k$-cycles for $k\geq 4$, we first note that the minimum number of induced $2k$-cycles that a vertex in $F_{n,2k}$ is contained in is $(n/k)^{k-1}+O(n^{k-2})$. In step 1 it would therefore be appropriate to consider $n$-vertex planar graphs in which every vertex is in at least $cn^{k-1}$ induced $2k$-cycles, where $c>0$ is a constant and $n$ is large. In Lemma~\ref{lem:longer_mainlem} below we show that such graphs contain structure analogous to that in Corollary~\ref{cor:good_cycle}, and similar to that found in the graph $F_{n,2k}$.

Given the analogues of steps 2 and 3 in the case of $2k$-cycles, an argument similar to that summarised in step 4 would yield Conjecture~\ref{conj:induced_even} for even cycles. Unfortunately, it does not seem to be straightforward to generalise steps 2 and 3. The main obstacle is that given a path $uvw$ in a graph $G$, for $k\geq 4$ it is possible that a vertex may appear in different `positions' in different induced $2k$-cycles containing the path $uvw$. For example there might exist one induced $2k$-cycle in which vertex $x$ is at distance $2$ from $u$ in the cycle and another in which it is at distance 3.

Nevertheless, below we give two results which roughly correspond to generalisations of Lemma~\ref{lem:empty} and Corollary~\ref{cor:good_cycle} to the case of $2k$-cycles for $k\geq 4$ in the hope that they might be of some use in a proof of Conjecture~\ref{conj:induced_even} in the case of even cycles, and at the very least to lend some credibility to Conjecture~\ref{conj:induced_even} in this case. For the lemmas to be relevant, there needs to exist for each $k\geq 4$ a constant $c>0$ such that if $n$ is sufficiently large relative to~$k$, then every vertex in an $n$-vertex planar graph containing $f_I(n,C_{2k})$ induced $2k$-cycles is in at least $cn^{k-1}$ induced $2k$-cycles. In light of the fact that $f_I(n,C_{2k})= (n/k)^k+o(n^k)$ (as noted in the introduction), this does not seem an unreasonable assumption.

\begin{lemma}\label{lem:empty_longer}
Let $k\geq 4$ be an integer, and let $c_0,c_1>0$ be constants. Let $n$ be large, and let $G$ be an $n$-vertex planar graph. Suppose that every vertex of $G$ is in at least $c_0n^{k-1}$ induced $2k$-cycles, and that there exist distinct vertices, $u$ and $w$, of $G$ with a common neighbourhood of size at least $c_1n$. Then there exists $v\in N(u)\cap N(w)$ such that every induced $2k$-cycle containing $v$ also contains $u$ and $w$.
\end{lemma}

The proof of Lemma~\ref{lem:empty_longer} is a straightforward adaptation of that of Lemma~\ref{lem:empty}, and we omit it. The only additional observation required for the proof is that in an $n$-vertex planar graph $G$ containing distinct vertices $u$ and $w$, the number of paths of length $i$ from $u$ to $w$ is at most $n(6n-12)^{i/2-1}$ if $i$ is even, and at most $(6n-12)^{(i-1)/2}$ if $i$ is odd.

\begin{lemma}\label{lem:longer_mainlem}
Let $k\geq 4$ be an integer and let $c>0$ be a constant. Then there exists a constant $c'>0$ such that if $n$ is large and $G$ is an $n$-vertex planar graph each of whose vertices is in at least $c n^{k-1}$ induced $2k$-cycles, then $G$ contains distinct vertices $x_0,x_2,x_4\dots,x_{2k-2}$ such that for each $i\in \{0,1,\dots,k-1\}$ we have $|N(x_{2i})\cap N(x_{2i+2})|\geq c'n$, where addition in the subscript is taken modulo $2k$.
\end{lemma}
\begin{proof}

Let $n$ and $G$ be as in the statement of the lemma, and let $v$ be a vertex of $G$ which has at most five neighbours with which it is in an induced $2k$-cycle in $G$. Such a vertex exists since $G$ is planar and thus contains a vertex of degree at most 5. Let $x_0$ and $x_{2k-2}$ be distinct neighbours of $v$ such that there are at least $cn^{k-1}/10$ induced $2k$-cycles in $G$ which contain the path $x_0vx_{2k-2}$.

\begin{claim}\label{claim:longer}
For each $i\in\{0,1,2,\dots,k-2\}$ and constant $\lambda>0$, there exists a constant $\lambda'=\lambda'(\lambda,i)>0$ such that if $x_1,x_2,x_3\dots,x_{2i}$ are distinct vertices such that the path $x_{2k-2}vx_0x_1\dots x_{2i}$ is contained in at least $\lambda n^{k-i-1}$ induced $2k$-cycles in $G$, then there exist vertices $x_{2i+2},x_{2i+4},x_{2i+6}\dots,x_{2k-4}$ all distinct from each other and from $x_0,x_1,\dots,x_{2i}$, and $x_{2k-2}$ such that for all $j\in\{i,i+1,\dots,k-2\}$ we have $|N(x_{2j})\cap N(x_{2j+2})|\geq \lambda' n$.
\end{claim}

\begin{proof}
We will prove the claim by reverse induction on $i$. If $i=k-2$, then let $\lambda>0$ and suppose that $x_1, x_2,\dots,x_{2k-4}$ are distinct vertices of $G$ such that the path $x_{2k-2}vx_0x_1\dots x_{2k-4}$ is contained in at least $\lambda n$ induced $2k$-cycles. This path has length $2k-2$, so $x_{2k-4}$ and $x_{2k-2}$ have at least $\lambda n$ common neighbours. Thus we can take $\lambda'(\lambda,k-2)=\lambda$.

Now let $i\in\{0,1,2,\dots,k-3\}$ and suppose that $\lambda'(\lambda,i+1)>0$ exists for all $\lambda>0$. Let $\lambda >0$, and let $x_1,x_2,\dots,x_{2i}$ be distinct vertices of $G$ such that the path $x_{2k-2}vx_0x_1\dots x_{2i}$ is contained in at least $\lambda n^{k-i-1}$ induced $2k$-cycles. Let $L$ be the set of ordered pairs $(a,b)$ of vertices of $G$ such that there is an induced $2k$-cycle in $G$ of the form $x_{2k-2}vx_0x_1\dots x_{2i}zaby_{2i+4}\dots y_{2k-3}$, and let $L'$ be the subset of $L$ consisting of the pairs which appear in at least $\lambda_0n^{k-i-2}$ cycles of this form, where $\lambda_0>0$ is a small constant to be determined later. Every induced $2k$-cycle containing the path $x_{2k-2}vx_0x_1\dots x_{2i}$ is of this form for some pair $(a,b)\in L$, and each pair appears in at most $n(6n)^{k-i-3}$ of these induced $2k$-cycles, because there are at most $n$ options for the vertex $z$, and at most $6n-12< 6n$ options for each of the edges $y_{2i+4}y_{2i+5},\dots,y_{2k-4}y_{2k-3}$. We also have $|L|\leq 6n-12<6n$, so if $l=|L'|$, then the total number of induced $2k$-cycles containing the path $x_{2k-2}vx_0x_1\dots x_{2i}$ is at most 
\[
ln(6n)^{k-i-3}+(6n-l)\lambda_0n^{k-i-2}.
\]
This expression is therefore at least $\lambda n^{k-i-1}$, so rearranging we obtain
\[
l\geq\frac{\lambda-6\lambda_0}{6^{k-i-3}-\lambda_0}n.
\]

Let $\lambda_1=(\lambda-6\lambda_0)/(6^{k-i-3}-\lambda_0)$ and choose $\lambda_0>0$ small enough that $\lambda_1>0$. For each pair $(a,b)\in L'$ there are at least $\lambda_0n^{k-i-2}$ induced $2k$-cycles in $G$ of the form $x_{2k-2}vx_0x_1\dots x_{2i}zaby_{2i+4}\dots y_{2k-3}$. There are at most $6n-12<6n$ options for each of the pairs $(y_{2i+4},y_{2i+5}),\dots,(y_{2k-4},y_{2k-3})$, and there are $k-i-3$ of these pairs, so the number of common neighbours of $x_{2i}$ and $a$ is at least $\lambda_0n/6^{k-i-3}$.

Let $\lambda_2=\lambda_0/6^{k-i-3}$ so that for each pair in $L'$, the vertex in the first position has degree at least $\lambda_2n$. Since the total number of edges in $G$ is at most $3n-6<3n$, there are at most $6/\lambda_2$ vertices which appear as the first vertex in a pair in $L'$. Let $\lambda_3=6/\lambda_2$. Since $|L'|\geq \lambda_1n$, there is a vertex, $x_{2i+2}$, of $G$ which is the first vertex in at least $\lambda_1n/\lambda_3$ pairs in $L'$. We have shown that $x_{2i}$ and $x_{2i+2}$ have at least $\lambda_2n$ common neighbours.

Moreover, there are at least $(\lambda_1\lambda_0/\lambda_3)n^{k-i-1}$ induced $2k$-cycles in $G$ of the form $x_{2k-2}vx_0x_1\dots x_{2i}zx_{2i+2}y_{2i+3}y_{2i+4}\dots y_{2k-3}$, so since $G$ contains only $n$ vertices, there exists $x_{2i+1}\in N(x_{2i})\cap N(x_{2i+2})$ such that there are at least $(\lambda_1\lambda_0/\lambda_3)n^{k-i-2}$ induced $2k$-cycles in $G$ containing the path $x_{2k-2}vx_0x_1\dots x_{2i}x_{2i+1}x_{2i+2}$. Applying the induction hypothesis, we see that we may take $\lambda'(\lambda,i)=\min(\lambda_2,\lambda'(\lambda_1\lambda_0/\lambda_3,i+1))$, and Claim~\ref{claim:longer} is proved.
\end{proof}

Applying the claim with $i=0$ and $\lambda=c/10$ shows that there exists a constant $c''>0$ and two distinct vertices $u$ and $w$ in $G$ which have a common neighbourhood of size at least $c''n$. By Lemma~\ref{lem:empty_longer}, there exists a vertex in $N(u)\cap N(w)$ such that every induced $2k$-cycle containing that vertex also contains $u$ and $w$. We may therefore take $v$ to be such a vertex, then take $x_0=u$ and $x_{2k-2}=w$ and apply the claim again with $i=0$ and $\lambda=c$ to complete the proof of Lemma~\ref{lem:longer_mainlem}.
\end{proof}

\section*{Acknowledgements}
Thank you to Alex Scott for advice on the drafting of this paper. Thank you also to an anonymous referee for their helpful comments and suggestions.

\end{document}